\documentclass[11pt]{amsart}
\usepackage{amsmath,amsthm, amssymb, amsfonts,amscd, epsfig}

\usepackage[mathscr]{eucal}

\theoremstyle{plain}

\newtheorem{thm}{Theorem}[section]
\newtheorem{lem}[thm]{Lemma}
\newtheorem{pro}[thm]{Proposition}
\newtheorem{cor}[thm]{Corollary}
\theoremstyle{definition}
\newtheorem{Def}[thm]{Definition}
\theoremstyle{remark}
\newtheorem{rem}[thm]{Remark}

\newcommand{\mc}[1]{\mathcal{#1}}

\newcommand{\mbf}[1]{\mathbf{#1}}

\newcommand{\mr}[1]{\mathrm{#1}}

\newcommand{\coh}[3]{ H ^{#1}  ( #2, \,  #3 )}

\newcommand{\lclcoh}[3]{H ^{#1}_{#2}  ( #3  )}

\newcommand{\Ass}{{\mathrm{Ass} }}

\newcommand{\Proj}{{\mathrm{Proj}\, }}

\newcommand{\Z}{{\mathbb Z}}
\newcommand{\R}{{\mathbb R}}

\newcommand{\proj}[1]{{\mathbf P}^{#1}}

\newcommand{\varlist}[3]{#1 _{#2}, \ldots ,#1 _{#3} }
\newcommand{\poly}[3]{#1 [\varlist{x}{0}{#2}, \varlist{y}{0}{#3}]}

\begin{document}

\title[Regularity]{Castelnuovo-Mumford regularity
in biprojective spaces}

\author{J. William Hoffman and Hao Hao Wang}
\address{Department of Mathematics\\
              Louisiana State University\\
              Baton Rouge, Louisiana 70803}

\email{hoffman@math.lsu.edu, wang\_h@math.lsu.edu}

\thanks{We would like to thank William Adkins
and David Cox for numerous discussions and
suggestions.}


\begin{abstract}
We define the concept of regularity for bigraded
modules over a bigraded polynomial ring. In this setting
we prove analogs of some of the classical results 
on $m$-regularity for graded modules over polynomial
algebras.  
\end{abstract}

\maketitle

\section{Introduction}
\label{S:intro}
In chapter 14 of \cite{DM} Mumford introduced the concept
of {\it regularity} for a coherent sheaf $\mc{F}$
on projective space $\proj{n}$: $\mc{F}$ is
$p$-regular if, for all $i \ge 1 $ we have
vanishing for the twists
\[
\coh{i}{\proj{n}}{\mc{F} (k)} = 0,\quad
\text{for all } k + i= p.
\]
This in turn implies the stronger condition of vanishing
for $k + i \ge p$. Regularity was investigated later
by several people, notably
Bayer and Mumford \cite{BM}, Bayer and Stillman \cite{BS},
Eisenbud and Goto \cite{EG}, and Ooishi \cite{O}. Let
$R = K[x_0 , ..., x_n]$ be the polynomial algebra in
$n+1$ variables over a field $K$,
graded in the usual way. If $M$ is a finitely generated
graded $R$-module, then the local cohomology groups
$\lclcoh{i}{\mbf{m}}{M}$ with respect to the ideal
$\mbf{m} = (x_0 , ..., x_n )$ are graded in a natural way
and we say that $M$ is $p$-regular if
\[
\lclcoh{i}{\mbf{m}}{M} _{k} = 0 \quad \text{ for all }
k+i \ge p+1.
\]
If $\mc{F}$ is the coherent sheaf on
$\proj{n}$ associated with $M$ in the usual way,
we have
\[
\lclcoh{i+1}{\mbf{m}}{M} _{k} =
\coh{i}{\proj{n}}{\mc{F} (k)}
\quad \text{for all } i \ge 1,
\]
which shows the compatibility of these definitions. An important
result in this theory is:
\begin{thm}
\label{T:BM1}
Suppose $K$ is a field and $I \subset R$ is a graded ideal.
Then $I$ is $p$-regular if and only if the minimal
free graded resolution of $I$ has the form
\[
\begin{CD}
0 @>>> \displaystyle{\bigoplus _{\alpha = 1}^{r_{s}}R e_{\alpha
,s}} @>>> \cdots @>>>\displaystyle{\bigoplus _{\alpha =
1}^{r_{1}}} R e_{\alpha ,1} @>>>  \displaystyle{\bigoplus _{\alpha
= 1}^{r_{0}}}R e_{\alpha ,0} @>>> I @>>> 0
\end{CD}
\]
where $\deg (e _{\alpha ,i}) \le p + i$ for all $i \ge 0$.
\end{thm}
The conditions of $p$-regularity can be derived quasi-axiomatically
from the following considerations. One seeks a condition in the form
of
\begin{equation}
\label{E:regdef}
\lclcoh{i}{\mbf{m}}{M}_{k} = 0 \text{ all } i \ge 0,
\text{ all }k \in C_i (p)
\Longrightarrow
R_s M_p = M_{p+s}\text{ all } s\ge 0,
\end{equation}
for certain regions $C_i (p) \subset \Z$.
One postulates:
\begin{itemize}
\item [1.] For each $i$, region $C_i (p)$ is independent
of the number $n+1$ of variables.
\item [2.] If $M$ is $p$-regular in the sense of
the left-hand side of
(\ref{E:regdef}),
then for a generic linear form $x \in R_1$,
$\bar{M} = M / xM$ is $p$-regular.
\end{itemize}First, when $n +1 = 0$, that is, we are considering
graded $K$-modules, since $\mbf{m} = (0)$, we have
$\lclcoh{0}{\mbf{m}}{M} = M$, and since $R_s = 0$ for
$s \ge 1$, property (\ref{E:regdef}) forces
$M_k =\lclcoh{0}{\mbf{m}}{M}_k = 0 $ for $k \ge p+1$
in this case. By principle 1., this must hold
for all $n$. Assuming that $\mbf{m}\notin \Ass(M)$ where $\Ass(M)$ 
denotes the associated primes for $M$, and $K$ is infinite, 
then $x$ may be chosen so that we have an exact sequence
\[
\begin{CD}
0@>>> xM = M(-1) @>>> M @>>> \bar{M}
@>>> 0
\end{CD}
\]
which gives rise to the long exact sequence in cohomology.
We have
\[
\begin{CD}
\lclcoh{0}{\mbf{m}}{M}_k @>>>
\lclcoh{0}{\bar{\mbf{m}}}{\bar{M}}_k
@>>>\lclcoh{1}{\mbf{m}}{M(-1)}_{k} =
\lclcoh{1}{\mbf{m}}{M}_{k-1}.
\end{CD}
\]
In order that we have
$\lclcoh{0}{\bar{\mbf{m}}}{\bar{M}}_k = 0$ for
$k \ge p+1$, as is demanded by principle 2., we must have
$\lclcoh{1}{\mbf{m}}{M}_{k} = 0$ for $k \ge p$. In a similar
way, we obtain the vanishing region for $\lclcoh{2}{\mbf{m}}{M}_{k}$
from that of $\lclcoh{1}{\mbf{m}}{M}_{k}$, etc., and we find
that they are exactly the conditions of $p$-regularity given.
Of course, one deduces property (\ref{E:regdef}) from the
definition of $p$-regularity, by induction on the number
of variables $n+1$, by a reversal of the above steps.
\par
The other essential feature of $p$-regularity is that
\begin{itemize}
\item[3.] $R$ is $0$-regular.
\end{itemize}
This follows from Serre's calculations of the cohomology
of the invertible sheaves $\mc{O}(k)$ on $\proj{n}$
(\cite{Serre1}),
as reinterpreted by Grothendieck in the language of local
cohomology (combine \cite[Prop. (2.1.5)]{aGIII} with
\cite[Exp. II, Prop. 5]{aG}).
\par
Our definition of regularity for bigraded modules follows this
pattern. Let $R = K[x, y] = K[x_0 , ..., x_m , y_0 , ..., y_n]$,
which is bigraded in the usual way. Let
$\mbf{m} = (xy) = (x_i y_j)$ be the irrelevant ideal.
We seek regions $C_{i}(p, p')\subset \Z ^2$ with the
property that

\begin{equation}
\label{E:regdef2}
\lclcoh{i}{\mbf{m}}{M}_{k, k'} = 0 \text{ all } i \ge 0,
\text{ all } (k, k')\in C_{i}(p, p')
\Longrightarrow
R_{s, s'} M_{p, p'} = M_{p+s, p'+s'}\text{ all } s\ge 0, s'\ge 0
\end{equation}
One postulates the analogs of 1. and 2. above. For 2.
we need regularity for both $M/xM$ and $M/yM$ for
generic $x \in R_{1, 0}$ and $y \in R_{0,1}$.
This leads to the regions called $Reg _{i-1}(p, p')$ (the
shift $i \to i-1$ is explained later). We are able to
prove analogs in this setting of the many of the classical
results for graded modules (see theorem (\ref{T:mainmodule})
and proposition (\ref{P:propweakreg})).  
Actually, we first do a separate
treatment for sheaves, the way Mumford did
(propositions (\ref{P:greater}) and (\ref{P:span})). However, in
attempting to generalize theorem (\ref{T:BM1}) to a
structure theorem for free resolutions for bigraded modules,
the conditions we have proposed are seen to be inadequate.
Therefore, we define a new concept of {\it strong} regularity
and prove that it does indeed give the structure theorem that
we want (theorem (\ref{T:BayerMum})). 
This involves vanishing conditions for all three of
$\lclcoh{\ast } {I}{M}$ for the ideals
$I = (x), (y), (x, y)$. The previous notion of regularity is now called
{\it weak} regularity. We show that strong regularity implies
weak regularity, and that $R$ itself is strongly $(0, 0)$-regular.
As far as we can determine, there is no simple vanishing condition
for  $\lclcoh{\ast } {(xy)}{M}$ alone that implies the structure
theorem that we want.
\par
In the last section we write down a free resolution that permits
computation of $\lclcoh{i}{\mbf{m}}{M}$. In a sequel to this work
applications and examples will be discussed. Also, it is clear that 
the methods in this paper may be extended to a multigraded module
over a multigraded polynomial algebra. This will also be addressed 
in a future work.  

\section{Regularity for coherent sheaves}
\label{S:sheafreg} First, we will give definition  and some
properties of regularity of a coherent sheaf similar to \cite[Ch.
14]{DM}. Let $K$ be a field, and $R=K[x_0, \cdots, x_m, y_0,
\cdots, y_n]$ be the polynomial ring, bigraded with variables $x$
having bidegree $(1, 0)$ and variables $y$ having bidegree $(0,
1)$.  We let
\[
\mbf{m} = R_{+} = \bigoplus _{a > 0, b > 0} R_{a, b},
\]
the {\it irrelevant ideal}. Some of the general theory of graded
and multigraded algebras used here can be found in \cite{GW1},
\cite{GW2}.
\par
Let $X=\proj{m} \times \proj{n}$, which when regarded as a scheme
is $\Proj (R)$, where by definition, this is the set of bigraded
prime ideals $\mbf{p}$ that do not contain the irrelevant ideal
$\mbf{m}$. There are projections $p_1$ and $p_2$ of $X$ onto its two
factors. If $\mc{F}_1$ is sheaf of $\mc{O}_{\proj{m}}$-modules,
and $\mc{F}_2$ is sheaf of $\mc{O}_{\proj{n}}$-modules, we denote
\[
\mc{F}_1 \boxtimes \mc{F}_2 = p_1 ^{\ast}\mc{F}_1 \otimes
 p_2 ^{\ast}\mc{F}_2 ,
\text{   an } \mc{O}_X \text{-module}.
\]
\par
As in the usual case of projective space there is a functor $M \to
\tilde {M}$ from bigraded $R$-modules to quasi-coherent sheaves on
$X$, and every quasi-coherent sheaf $\mc{F}$ arises this way, in a
nonunique fashion. In fact, if
\[
M = \bigoplus _{(a, b)\in \Z ^{2}}\coh{0}{X}{\mc{F}(a, b)}
\]
then $\mc{F} \cong \tilde{M}$. Here, for any sheaf of
$\mc{O}_{X}$- modules $\mc{F}$, we denote
\[
\mc{F}(a, b) = \mc{F}\otimes \mc{O}_{X} (a, b)
\]
where $\mc{O}_{X} (a, b)= \mc{O}_{\proj{m}}(a)
\boxtimes\mc{O}_{\proj{n}}(b) $ is the invertible sheaf associated
to the graded $R$-module $R(a, b)$. Recall that if $M$ is any
graded $R$-module, $M(a, b)$ is the graded module with degrees
shifted via $M(a, b)_{d, e} = M_{d+a, e+b}$. If $Z$ is a scheme,
tensor products involving $\mc{O}_{Z} $-modules will always be
relative to $\mc{O}_{Z}$ unless otherwise stated
\par
When $m\geq 1$, and $n\geq 1$, the Picard group $\mr{Pic} (X)$ is
isomorphic with $\Z ^2$ with $(a, b)$ corresponding to $\mc{O}_{X}
(a, b)$. Interpreting the Picard group as the group of
divisor-classes, $\mc{O}_{X} (a, b)$ corresponds to the divisor $a
L_1 + b L_2$, where $L_1 = H_1 \times \proj{n}$, $H_1 \subset
\proj{m}$ being any hyperplane, and $L_2 =  \proj{m}\times H_2$,
$H_2 \subset \proj{n}$ being any hyperplane.
\par
Note the special case: if $m$ or $n$ is $0$, the biprojective
space reduces to a projective space.  Except in the case where
both are $0$, the Picard group $\mr{Pic}(X)$ is isomorphic with
$\Z$. If both are $0$, the space reduces to a point, and its
Picard group is trivial. Even in these degenerate cases we still
use notations such as $\mc{F}(a, b)$, where one or other twisting
by $a$ or $b$ might be trivial.

\begin{Def}
\label{D:stair}
For  each integer $i > 0$, let
\begin{eqnarray*}
St_i &=& \{(r,s)\in \Z^2: r+s=-i-1, \ r<0, \ s<0\}\\
&= &\{(-i, -1),\, (-i+1, -2),\, \ldots ,
         \, (-2, -i +1),\, (-1, -i)
       \},
\end{eqnarray*}
for $i \leq 0$, let
\begin{eqnarray*}
St_i &=& \{(r,s)\in \Z^2: r+s=-i, \ r\geq 0, \ s \geq 0\}\\
&= &\{(-i,0),\, (-i-1, 1),\, \ldots ,
         \, (1, -i-1),\, (0,-i)
       \}.
\end{eqnarray*}
For each $(p,p')\in \Z^2$ let $St_i (p, p') = (p, p') + St _i$.
\newline
For $i \geq 0$, let $Reg_i (p, p' ) = \Z _{+} ^2 + St_i (p, p')$
where $\Z_+= \{ n \in \Z : n \ge 0 \}$. \newline
For $i=-1$, let $Reg_{-1}=\Z_+^2+(p+1,p'+1)$.\newline
Let  $Reg' _{-1} (p, p' ) = (p +1, p') + \Z _{+} ^2$.\newline
Let $Reg'' _{-1} (p, p' ) = (p, p' + 1) + \Z _{+} ^2$.\newline
For $i \geq 0$, define
$DReg_i(p,p')=\Z^2_-+St_{-i}(p,p')$  where $\Z_-=\{n \in \Z: n
\leq 0\}$.\newline
Note that, for all $i\ge -1$,
\[
Reg_i (p, p') \ = \ \Z ^2 \cap \{(x, y)\in \R ^2 \mid
         x \ge p-i , \ y \ge p'-i, \
         x +y \ge p+p'-i-1   \}
\]
and, for all $i\ge 0$,
\[
DReg_{i}(p, p') = -Reg_{i+1}(-p+1, -p'+1)
\]
\end{Def}

\begin{rem}\label{R:streg}
For $i \geq 0$, and for all $p,p'$, we have
\begin{itemize}
\item[1.] $(k,k')\in St_i(p,p') \Rightarrow (k-1,k'), (k,k'-1) \in St_{i+1}(p,p')$.
\item[2.] $St_i(p,p')\in Reg_i(p,p')$.
\item[3.] $(k,k')\in Reg_i(p,p') \Rightarrow (k-1,k'), (k,k'-1)
\in Reg_{i+1}(p,p')$.
\item[4.] $Reg_i(q,q') \subset Reg_i(p,p')$, if $q \geq p, q'\geq p'$.
\item[5.] $(k, k') \in Reg' _{-1} (p, p') \Longrightarrow (k-1, k')
\in Reg_{0} (p, p')$.
\item[6.]  $(k, k') \in Reg''_{-1} (p, p') \Longrightarrow (k, k'-1)
\in Reg_{0} (p, p')$.
\end{itemize}
\end{rem}

Here is a picture:
\begin{figure}[h]
\label{F:1}
\input{HReg.pstex_t}
\caption{$Reg_{i}(p, p')$}
\end{figure}

\begin{figure}[h]
\label{F:2}
\input{HDReg.pstex_t}
\caption{$DReg_i(p,p')$}
\end{figure}

Using these notations, we make the following definition.
\begin{Def}
\label{D:sheafreg} Let $\mc{F}$ be a coherent sheaf on $X$. We
will say that $\mc{F}$ is \textit{$(p,p')$-regular} if, for all $i
\ge 1$,
\[
H^i(X, \mc{F}(k,k'))= 0
\]
whenever $(k,k')\in St_i(p,p')$.
\end{Def}

\begin{rem} \label{R:Reduced}
If $n=0$, $\proj{m}\times \proj{0} \cong \proj{m}$, so every
coherent sheaf on $\proj{m}\times \proj{0}$ is naturally
identified with a sheaf on $\proj{m}$. The sheaf $\mc{F}(p,p')$ is
independent of $p'$. Under this identification, $\mc{F}$ is
$(p,p')$-regular on $\proj{m}\times \proj{0}$ in the sense of Definition
\ref{D:sheafreg}, if and only if $\mc{F}$ is $p$-regular on
$\proj{m}$ in the sense of Mumford.
\end{rem}
\begin{proof}
First, we  will show that $(p,p')$-regular implies $p$-regular.
\par
 In this case, $\mc{F}(k,k')\cong \mc{F}(k)$.  $\mc{F}$ is
$(p,p')$-regular means that for all $i \geq 1$,
\[
H^i(\proj{m}\times \proj{0}, \mc{F}(k,k'))=H^i(\proj{m},
\mc{F}(k))=0,
\]
where $p-i \leq k \leq p-1$. Since $k+i \geq p$, according to
\cite[p. 100]{DM}, $\mc{F}$ is $p$-regular.
\par
Second, we will show that $\mc{F}$ is $p$-regular implies
$(p,p')$-regular.
\par
If $\mc{F}$ is $p$-regular, then $H^i(\proj{m}, \mc{F}(k))=0$
whenever $k+i \geq p$,  this implies
$$H^i(\proj{m}\times \proj{0}, \mc{F}(k,k'))=0$$
for any $k' \in \Z$. In particular, $H^i(\proj{m}\times \proj{0},
\mc{F}(k,k'))=0$ for all $(k,k') \in St_i(p,p')$.  Therefore,
$\mc{F}$ is $(p,p')$-regular.
\end{proof}

\begin{pro}
\label{P:zerozero} $\mc{O}_{X}$ is $(0, 0)$-regular.
\end{pro}
\begin{proof}
  If $m$ or $n=0$, by the previous
remark, $\mc{O}_X$ is $(0,0)$-regular $\Leftrightarrow$ $\mc{O}_X$
is $0$-regular. But $\mc{O}_{\proj{m}}$ is $0$-regular since
\begin{eqnarray} \label{serre1}
 \coh{a}{\proj{m}}{\mc{O}_{\proj{m}}(k)} &=& 0,
\text{ if } a\geq 1 \text{ and } a + k \geq 0\\
\label{serre2} \coh{0}{\proj{m}}{\mc{O}_{\proj{m}}(k)} &=& 0,
\text{ if } k \leq -1.
\end{eqnarray}
These formulas are a consequence of Serre's results on the
cohomology of projective space. \cite{H}

 If $m$ and $n\ge 1$, we can apply the K\"unneth formula \cite{SW},
\[
\coh{i}{X}{\mc{O}_{\proj{m}}(k) \boxtimes\mc{O}_{\proj{n}}(k')} =
\bigoplus _{a + b = i}
\coh{a}{\proj{m}}{\mc{O}_{\proj{m}}(k)}\otimes
                         \coh{b}{\proj{n}}{\mc{O}_{\proj{n}}(k')}.
\]
We will show that $ \coh{a}{\proj{m}}{\mc{O}_{\proj{m}}(k)}=0$ or
$\coh{b}{\proj{n}}{\mc{O}_{\proj{n}}(k')}=0$ whenever $a+b=i$ and
$(k,k')\in St_i(0,0)$. If $(k,k') \in St_i(0,0)$, then  $k=-i+l$
and $k'=-1-l$ where $0 \leq l \leq i-1$. If $a=0$ or $b=0$, we are
done by Equation \eqref{serre2}, since $k, k' <0$. If $a
>0$, and $b>0$, we only need to show $a-i+l \geq 0$ or $b-1-l \geq
0$.   Suppose both $a-i+l \leq -1 $ and $b-1-l \leq -1 $.  Since
$a+b=i$, $$-1=(a-i+l)+(b-1-l)\le -2.$$ This contradiction shows
that either $a-i+1\ge 0$ or $b-1-l\ge 0$, and the proof is
completed by Equation \eqref{serre2}.
\end{proof}

\begin{lem}
\label{L:restrict} Assume that $K$ is infinite, and that $m\ge 1$.
Let $\mc{F}$ be a
 coherent sheaf on $X$. Let $L_1$ be a hyperplane
defined by $\sum_{i=0}^m a_i x_i=0$, and let $\mc{F}_{L_1}=
\mc{F}\otimes \mc{O}_{L_{1}}$ denote the sheaf $\mc{F}$ restricted
to $L_1$. If $\mc{F}$ is $(p,p')$-regular, then $\mc{F}_{L_1}$ is
$(p,p')$-regular for a generic $L_1$. The similar statement is
true for hyperplanes $L_2$ defined by a
 form $\sum_{i=0}^n b_i y_i=0$ assuming $n \geq 1$.
\end{lem}

\begin{proof}
Given $\mc{F}$, choose a hyperplane $L_1$, where $L_1$ is defined
by an equation of the form $f = \sum^{m}_{i=0}a_i x_i=0$, such
that $L_1$ does not contain any of points of the finite set of
associated primes $A (\mc{F})$ (for the definition of this, see
\cite[p.40]{DM}). Note that this is possible:  $A (\mc{F})$ is
finite, and because $K$ is infinite, we can find a linear form
missing the $p_1 $-projections of the associated primes. Tensor
the exact sequence
\[
\begin{CD} 0   @>>>  \mc{O}_X(-1,0)  @> f >>  \mc{O}_X \rightarrow
\mc{O}_{L_1}  @>>>  0
\end{CD}
\]
with $\mc{F}(k,k')$.  For all $x \in X$,  multiplication by $f$ is
 injective in $\mc{F}_{x}$, since by construction, $f$ is a
 unit at all associated primes of $\mc{F}_x$. Therefore the
 resulting sequence is exact:
 \begin{equation}
\label{E:exact1}
\begin{CD}
 0  @>>> \mc{F}(k-1,k') @>f>>
 \mc{F}(k,k') @>>>
\mc{F}\otimes \mc{O}_{L_1}(k,k')=  \mc{F}_{L_1}(k,k') @>>> 0
\end{CD}
\end{equation}
This gives an exact cohomology sequence:
\[
\begin{CD}
\cdots @>>> H^i(\mc{F}(k,k')) @>>> H^i(\mc{F}_{L_1}(k, k')) @>>>
H^{i+1}(\mc{F}(k-1,k'))@>>> \cdots
\end{CD}
\]

If  $(k, k') \in St _i (p, p')$, then $(k-1,k') \in
St_{i+1}(p,p')$ by \ref{R:streg}, and the first and the last
groups vanish when $i \geq 1$, since we are assuming that $\mc{F}$
is $(p, p')$-regular. This forces the second group to vanish, thus
proving that $\mc{F}_{L_1}$ is $(p,p')$-regular.
\end{proof}

\begin{pro}\label{P:greater}
 If $\mc{F}$ is a $(p,p')$-regular coherent sheaf on
$X=\proj{m} \times \proj{n}$, then for all $i \geq 1$,
\begin{equation} \label{E:RegReg}
H^i(X, \mc{F}(k,k'))=0
\end{equation}
whenever  $(k, k') \in Reg_i (p, p')$. That is, $\mc{F}$ is
$(q,q')$-regular for $q \geq p,\  q' \geq p'$.
\end{pro}

\begin{proof}
We will prove \eqref{E:RegReg} by double induction on $(m,n)$.  If
$m=0$ or $n=0$, by Remark \ref{R:Reduced} $(p,p')$-regularity
reduces to ordinary $p$-regularity or $p'$-regularity for
projective space, and \eqref{E:RegReg} is true by Mumford's result
\cite{DM}. So assume $m \ge 1$ and $n \ge 1$. Every element of
$Reg_i (p, p')$ is of the form $(k+ r, k'+s)$ for some $(k, k')\in
St_i (p, p')$, and $(r, s)\geq (0, 0)$. Now we will do double
induction on the pair $(r, s)$.  The case $(r, s)= (0, 0)$ is true
by assumption of $(p, p')$-regularity for $\mc{F}$. Choose a
hyperplane $L_1$ as in Lemma \ref{L:restrict} such that
$\mc{F}_{L_1}$ is $(p,p')$-regular.  Consider the cohomology exact
sequence attached to \eqref{E:exact1} with $(k, k')$ replaced by
$(k+r+1, k'+s)$:
\[
\begin{CD}
 H^i(\mc{F}(k+r,k'+s)) @>>>
H^i(\mc{F}(k+r+1, k'+s)) @>>> H^{i}(\mc{F}_{L_{1}}(k+r+1,k'+s))
\end{CD}
\]
Since $\mc{F}_{L_1}$ is $(p,p')$-regular, and since $L_1$ is a
biprojective space of lower dimension, the induction hypothesis
says that the right-hand term is $0$. The left-hand side also
vanishes, by induction hypothesis on $(r, s)$. Hence the middle
term vanishes, as required. A symmetric argument shows vanishing
for $(k+r, k'+s+1)$.
\end{proof}

\begin{pro}\label{P:span}
If $\mc{F}$ is a $(p,p')$-regular coherent sheaf on $X$, then
\item $H^0(X, \mc{F}(k,k'))$ is
spanned by
\[ H^0(X, \mc{F}(k-1,k'))
 \otimes H^0(X, \mc{O}(1,0)), \]
if $k > p, k' \geq p'$ ; and  it is spanned by
\[  H^0(X, \mc{F}(k,k'-1))
 \otimes H^0(X, \mc{O}(0,1)), \]
if $ k\geq p, k'>p'$.
\end{pro}

\begin{proof}
We use induction on $\dim(X)$: for $\dim(X)=0$, the result is
true. By Lemma \ref{L:restrict}, we know that $\mc{F}_{L_1}$ is
$(p,p')$-regular. Consider the following diagram:
\[
\begin{CD}
      @.    H^0(\mc{F}(k-1,k'))\otimes H^0(\mc{O}_X(1,0))
@>{\sigma}>> H^0(\mc{F}_{L_1}(k-1, k'))\otimes H^0(\mc{O}_{L_1}(1,0)) \\
      @.    @VV{\mu}V           @VV{\tau}V            \\
H^0(\mc{F}(k-1,k'))
@>\alpha>>H^0(\mc{F}(k,k'))@>{\nu}>>H^0(\mc{F}_{L_1}(k,k'))
\end{CD}
\] If $k > p$ and $k' \geq p'$,  $\sigma$ is
surjective because $\mc{F}$ is $(p,p')$-regular, and thus
$H^1(\mc{F}(k-2,k'))=0$. $\tau$ is surjective by induction
hypothesis. $\nu$ is also surjective, since $H^1(\mc{F}(k-1,
k))=0$.
\par
Let $t \in H^0(\mc{F}(k,k'))$, we have $\nu(t)=\tau(s)=\tau
\sigma(s')$ for some
\[
s \in H^0(\mc{F}_{L_1}(k-1, k'))\otimes H^0(\mc{O}_{L_1}(1,0)),
\text{   and } s' \in H^0(\mc{F}(k-1,k'))\otimes
H^0(\mc{O}_Q(1,0)).
\]
We have $\nu(\mu(s'))=\tau(\sigma(s'))=\nu(t)$, and $t-\mu(s')\in
ker(\nu)$. Since the last row of the diagram is exact in the
middle, so we have $t' \in H^0(\mc{F}(k-1,k'))$ such that
$\alpha(t')=t-\mu(s')$.  This says that $H^0(\mc{F}(k,k'))$ is
spanned by the image of $\mu$ and the image of $\alpha$. But the
image of $\alpha$ is in $H^0( \mc{F}(k-1,k')) \otimes
H^0(\mc{O}(1,0))$, because map $\alpha$ is the multiplication by
$f$, and $f \in H^0(\mc{O}(1,0))$.  This means that
$H^0(\mc{F}(k,k'))$ is spanned by
\[ H^0(\mc{F}(k-1,k')) \otimes H^0( \mc{O}(1,0))\]
 By symmetry, we can show that if $k \geq p, k' >p'$,
 $H^0(\mc{F}(k,k'))$ is spanned by
 \[ H^0( \mc{F}(k,k'-1)) \otimes H^0( \mc{O}(0,1)).\]
 \end{proof}

\section{Weak regularity for bigraded modules}
\label{S:modreg1}
We will give the definition and some properties
of regularity for a bigraded module similar to  \cite{O} and
\cite{JK}. Let $A$ be a noetherian ring, and let now
$R=\oplus_{a,b \geq 0} R_{a,b}$ be any bigraded ring over $A$,
with $R_{0,0}=A$. We assume that it is finitely generated by
homogeneous elements of bidegrees $(1, 0)$ and $(0, 1)$. Such a
ring will be called a bihomogeneous $A$-algebra. Previously we
considered only the case of a polynomial ring in two sets of
variables over a field.  Let $\mbf{m}=R_+=\oplus_{a>0,
b>0}R_{a,b}$ be the irrelevant ideal; it is a bigraded $R$-module.
There is a scheme $X= \Proj (R)$, whose points are the
bihomogeneous prime ideals $\mbf{p}$ of $R$ that do not contain
the irrelevant ideal. We also have a functor $M \to \tilde{M}$
from bigraded modules to quasicoherent $\mc{O}_X$-modules with
similar properties to those discussed in section \ref{S:sheafreg}.
Let $\mc{F}$ be a quasicoherent $\mc{O}_X$-module.  If we set
$M=\oplus_{a,b \in \Z}H^0(X, \mc{F}(a,b))$, then we have
$\mc{F}=\tilde{M}$.
\par
If $R$ is a bigraded $A$-algebra, then it defines a graded
$A$-algebra
\[
R^{\sharp}_n = \bigoplus _{i+j=n}  R_{i, j}
\]
and similarly we have a graded $R^{\sharp}$-module $M^{\sharp}$
associated to a bigraded $R$-module $M$.
\par
Let $M=\oplus_{a,b\in \Z}M_{a,b}$ be a bigraded $R$-module. The
local cohomology groups $H^i_{\mbf{m}}(M)$ are bigraded
$R$-modules, and let $\lclcoh{i}{\mbf{m}}{M}_{a, b}$ denote the
$(a,b)$ part. The general theory of local cohomology is found
in \cite{aG}. Note that, if $J \subset A$ is an ideal in a ring,
and $V(J)\subset C = \text{Spec}(A)$ is the corresponding
closed subset, then
\[
\lclcoh{\ast}{J}{M} = \lclcoh{\ast}{V(J)}{C,\, \tilde{M}}
\]
where $\tilde{M}$ is the quasi-coherent sheaf on $C$
associated with the $A$-module $M$.
\par
We have
\[
\lclcoh{i}{\mbf{m}}{M}^{\sharp} =
\lclcoh{i}{\mbf{m}^{\sharp}}{M^{\sharp}},
\text{ ie., } \lclcoh{i}{\mbf{m}^{\sharp}}{M^{\sharp}}_{n} =
\bigoplus _{k + k' = n}
\lclcoh{i}{\mbf{m}}{M}_{k, k'}
\]
Generally we omit the $\sharp$ from $\mbf{m}$, as it is clear
in context that we are referring to the graded, as opposed to
the bigraded, structure.
\par
We recall the following fact: Let $R$ be any ring, $I\subset R$ an
ideal and $M$ an $R$-module. If $\mr{Supp}(M)\subset V(I)$ then
\[
\lclcoh{0}{I}{M} = M, \  \text{and} \ \lclcoh{i}{I}{M} = 0 \text{
for }i \ge 1.
\]
Also, if $R$ is Noetherian and $M$ is finitely generated,
\[
\mr{Ass}(M) \subset \mr{Supp}(M),
\]
and both have the same minimal elements. $\mr{Ass}(M)$ is finite.
\par
We allow the case where $R_{a, b} = 0$ for all $a > 0$, or $R_{a,
b} = 0$ all $ b >0$. For then $\mbf{m} = 0$, and thus for all
$R$-modules $M$,
\[
\lclcoh{0}{\mbf{m}}{M} = M, \  \text{and} \ \lclcoh{i}{\mbf{m}}{M}
= 0 \text{ for }i \ge 1.
\]
since $V(\mbf{m}) = \mr{Spec}(R)$, so $\mr{Supp}(M)\subset
V(\mbf{m})$ always holds. This extreme case plays an
important role in the proofs of the main theorems about
regularity, which are by induction on the number of variables.

\begin{Def}
\label{D:weakmodreg}
We say that a bigraded $R$-module $M$ is \textit{weakly
$(p,p')$-regular}, if for all $i \geq 0$,
\[
H^i_{\mbf{m}}(M)_{k,k'}=0  \text{  for all  } (k,k') \in
Reg_{i-1}(p, p')
\]
\end{Def}
The connection with the previous concept of regularity
for coherent sheaves is established by the following:

\begin{pro}
\label{P:exact}
(see \cite{HY}) Let $X = \mr{Proj}(R)$.
For any finitely generated bigraded $R$-module $M$
we have an exact sequence of bigraded $R$-modules
\[
\label{E:local1}
\begin{CD}
0 @>>> H^0_{\mbf{m}}(M) @>>> M @>>>
\displaystyle{\bigoplus _{(a,b) \in \Z ^2}H^0(X,\,  \mc{M}(a,b))} @>>>
H^1_{\mbf{m}}(M) @>>>0
\end{CD}
\]
and an isomorphism of bigraded $R$-modules
\[
H^{i+1}_{\mbf{m}}(M) = \bigoplus_{(a,b) \in \Z ^2} H^i (X, \, \mc{M}(a,b)),
\  \ \forall i \geq 1
\]
\end{pro}

\begin{cor}
\label{C:exact}
Let $\tilde{M}$ be the sheaf on $X$ associated to the
bigraded $R$-module $M$,
if $M$ is weakly $(p,p')$-regular, then $\tilde{M}$ is
$(p,p')$-regular in the sense of definition
\ref{D:sheafreg}. This explains the shift in index
from $i$ to $i-1$ in the definition of weak regularity
for modules.
\end{cor}

The main result for weak regularity is the following:

\begin{thm}
\label{T:mainmodule}
Let $R$ be bihomogeneous $A$-algebra,
$M$ a finitely generated bigraded
$R$-module. Fix $(p, p')$.
\begin{itemize}
\item[1.] Suppose that
$H^i_{\mbf{m}}(M)_{k,k'}=0$ for all $i \geq 1$
and all $(k,k') \in St_{i-1} (p, p')$, then
\[
H^i_{\mbf{m}}(M)_{k,k'}=0
\text{ for all } i \geq 1 \text{ and all }
(k,k') \in Reg_{i-1}(p, p')
\]
\item[2.]
Moreover,
\begin{itemize}
\item[a.] if $H^0_{\mbf{m}}(M)_{k,k'}=0$ for
$(k,k')\in Reg'_{-1} (p, p')$, then we have
$R_{d,0}M_{k,k'}=M_{d+k,k'}$ for every $d\ge 0,
k \geq p, k'\geq p'$;
\item[b.]
if $H^0_{\mbf{m}}(M)_{k,k'}=0$ for
$(k,k')\in Reg''_{-1} (p, p')$, then we have
 $R_{0,d'}M_{k,k'}=M_{k,k'+d'}$ for every
$d'\ge 0, k \geq p, k'\geq p'$.
\end{itemize}
\item[3.]
 if $M$ is weakly $(p, p')$-regular, and if it satisfies
$H^0_{\mbf{m}}(M)_{k,k'}=0$ for $(k,k') \in Reg_{-1}^{'}(p,p')\cup
Reg_{-1}^{''}(p,p')$, then
 $R_{d,d'}M_{k,k'}=M_{k+ d,k'+d'}$ for all
$d, d'\ge 0$,\newline
 $k \geq p, k'\geq p'$.
\end{itemize}
\end{thm}

\begin{proof}
First, by the same argument as in \cite[Theorem 2]{O}, we
may reduce to the case where $A$ is a local ring with
infinite residue field, and assume that
$R=A[x_0, \cdots, x_m, y_0, \cdots, y_n]$, with
irrelevant ideal $\mbf{m}$ generated by the
$x_i y_j$.
We will prove the claim by induction
on $(m, n)$.  If either $m = -1$ or $n=-1$ (ie., either
$x$ or $y$ variables are missing), or if
\[
\Ass_+(M)=\{ \mbf{p} \in \Ass(M):
\mbf{p} \nsupseteq \mbf{m}\}=\emptyset
\]
the claim is true: in the first case the irrelevant
ideal $\mbf{m} = 0$, so that the remark
before the statement of proposition \ref{P:exact} applies;
in the second case, we have $\mr{Supp}(M) \subset V(\mbf{m})$.
In either case,
$H^0_{\mbf{m}}(M)=M$ and $H^i_{\mbf{m}}(M)=0$ for every $i\geq 1$.
\par
Suppose that both $m\ge 0$ and $n \ge 0$,
and $\Ass_+(M)=\{ \mbf{p}_1, \cdots, \mbf{p}_r\}$.
By our assumptions, $A$ is a um-ring in the
terminology of \cite{QB}, and by theorem 2.3 of that paper we
conclude that if we had an equality of $A$-modules
\[
R_{1, 0}
=
\text{max}(A)R_{1, 0}\cup ( \mbf{p}_1\cap R_{1, 0})
\cup \cdots \cup (\mbf{p}_r\cap R_{1, 0})
\]
then $R_{1, 0}$ would have to be equal to one of the
terms in the union. It clearly is not the first term.
If, say  $R_{1, 0} = \mbf{p}_1\cap R_{1, 0}$
we would have
\[
\mbf{m} \subset (x_0 , ..., x_m) = R\cdot R_{1, 0}
\subset \mbf{p}_1
\]
which is contrary to the fact that $\mbf{p}_1 $ does
not contain $\mbf{m}$.
Thus we can find an element
\[
x \in R_{1, 0}
-
\text{max}(A)R_{1, 0}\cup ( \mbf{p_1}\cap R_{1, 0})
\cup \cdots \cup (\mbf{p_r}\cap R_{1, 0})
\]
which we can take as part of a free basis of $R_{1,0}$, and
by change of coordinate, we may assume that $x=x_m$.
\par
(1.)
Consider the following exact sequence:
\[
\begin{CD} 0 @>>> M_1 @>>>M \  @ > {x} >>
xM(1,0) @>>> 0. \end{CD}
\]
This implies:
\begin{equation}
\label{E:modseq}
\begin{CD}
H^i_{\mbf{m}}(M_1) @>>> H^i_{\mbf{m}}(M) @>>>
H^i_{\mbf{m}}(xM(1,0)) @>>> H^{i+1}_{\mbf{m}}(M_1).
\end{CD}
\end{equation}
Since $x$ was chosen not to belong to any of the
$\mbf{p}_i$, $\text{Supp}(M_1) \subset V(\mbf{m})$,
and so by the remarks above, the first and last terms above
vanish when $i\ge 1$, and so $H^i_{\mbf{m}}(M) \cong
H^i_{\mbf{m}}(xM(1,0))$ for every $i \geq 1$.  Set
$\bar{R}=R/xR=A[x_0, \cdots, x_{m-1},y_0, \cdots, y_n]$,
$\bar{\mbf{m}}=\bar{R}_+$ and $\bar{M}=M/xM$.
From
\[ \begin{CD}
0 @>>> xM @>>> M @>>> \bar{M} @>>>
0 ,\end{CD}
\]
we have the exact sequence:
\begin{equation}
\label{E:modseq2}
\begin{CD}
H^i_{\mbf{m}}(M)_{k,k'} @>>>
H^i_{\bar{\mbf{m}}}(\bar{M})_{k,k'} @>>>
H^{i+1}_{\mbf{m}}(xM)_{k,k'} =
H^{i+1}_{\mbf{m}}(M)_{k-1,k'}.\end{CD}
\end{equation}
If $(k,k') \in St_{i-1}(p, p')$, then the first term is $0$, by
our assumption on $M$.
\par
Now assume that $i \ge 2$. Then, $(k-1, k') \in St_{i}(p, p')$ by
remark (\ref{R:streg}), and so the last term above is $0$, also by
our assumption on $M$, so that
$H^i_{\bar{\mbf{m}}}(\bar{M})_{k,k'}=0$. By induction hypothesis
$H^i_{\bar{\mbf{m}}}(\bar{M})_{k,k'}=0$ for every $i \geq 1$ and
$(k,k')\in Reg_{i-1}(p, p')$. If $i \geq 1$ and $(k,k') \in St
_{i-1} (p + 1, p')$, then in the exact sequence
\begin{equation}
\label{E:modseq3}
\begin{CD}
H^i_{\mbf{m}}(M) _{k-1, k'}=H^{i}_{\mbf{m}}(xM)_{k,k'}
 @>>> H^{i}_{\mbf{m}}(M)_{k,k'} @>>>
H^i_{\bar{\mbf{m}}}(\bar{M})_{k,k'}\end{CD}
\end{equation}
the first and last terms are $0$ because
$(k-1, k')\in St_{i-1} (p, p')$ and $(k, k')
\in Reg _{i-1} (p, p')$,
so $H^{i}_{\mbf{m}}(M)_{k,k'}=0$ when
$i\ge 1$ and $(k,k')\in St_{i-1}(p+1,p')$.
Repeating the argument we get
 $H^{i}_{\mbf{m}}(M)_{k,k'}=0$ when
$i\ge 1$ and $(k,k')\in St_{i-1}(p+d,p')$ for every
$d \ge 0$, and by symmetry, arguing with
a $y \in R_{0, 1}$, we get $H^{i}_{\mbf{m}}(M)_{k,k'}=0$ when
$i\ge 1$ and $(k,k')\in St_{i-1}(p+d,p'+ d')$ for every
$d, d' \ge 0$, which is the first claim for $i\ge 2$.
\par
When $i =1$, the only changes to make in the argument
are the following. If  $(k,k')\in St_{0}(p,p')$, then
$(k-1,k')\in Reg_{1}(p,p')$, by formula (\ref{R:streg}).
But then  $H^{2}_{\mbf{m}}(M)_{k-1,k'}=0$ has been
established by the argument in the previous paragraph.
Also, when  $(k,k')\in St_{0}(p+1,p')$, we have
$(k-1, k')\in St _{0} (p, p')$ and $(k, k') \in Reg _{0} (p, p')$,
so that the first and last terms in the sequence
(\ref{E:modseq3})
vanish when $i=1$, too.

\par
(2a.) Let $\Ass_+(M)=\{ \mbf{p} \in \Ass(M):
\mbf{p} \nsupseteq \mbf{m}\}$. Suppose $m, n \ge0$ and
$\Ass_+(M)=\{
\mbf{p}_1, \cdots, \mbf{p}_r\}$. As before, we change
coordinates so that $x=x_m \notin \mbf{p}_i$, for any $i$.  Set
$\bar{R}=R/xR=A[x_0, \cdots, x_{m-1},y_0, \cdots, y_n]$,
$\bar{\mbf{m}}=\bar{R}_+$ and $\bar{M}=M/xM$. We claim
that the induction hypothesis can be applied to $\bar{M}$.
First, by the argument proving (1.), we saw that
\[
 H^i_{\bar{\mbf{m}}}(\bar{M})_{k,k'}= 0  \text{ for } i\ge 1
\text{ and } (k, k')\in Reg_{i-1}(p, p').
\]
From the
sequence (\ref{E:modseq2}) above with $i=0$, we see that
$H^0_{\bar{\mbf{m}}}(\bar{M})_{k,k'}=0$ for every
$(k, k') \in Reg'_{-1}(p, p')$,
because the extreme terms vanish: the left-hand
one because of our assumption on $M$, the right-hand one
because $(k-1, k')\in Reg_{0}(p, p')$ by remark
(\ref{R:streg}) and vanishing of this term has been
established above.
Thus by induction hypothesis applies to $\bar{M}$,
and we have
$\bar{R}_{d,0}\bar{M}_{k,k'} = \bar{M}_{d+k,k'}$,
which implies $R_{d,0}M_{k,k'}+xM_{d+k-1,k'}=M_{d+k,k'}$.
Reasoning by
induction on $d\ge 1$, we assume that
$M_{d+k-1,k'}=R_{d-1,0}M_{k,k'}$ has been established,
the case $d=1$ being trivial. Then
\[
M_{d+k,k'}=R_{d,0}M_{k,k'}+xM_{d+k-1,k'}=
R_{d,0}M_{k,k'} + xR_{d-1,0}M_{k,k'}=R_{d,0}M_{k,k'}.
\]
This proves our claim.
By symmetry, arguing with a $y_n$, we get the assertion
$M_{k,d'+k'}=R_{0,d'}M_{k,k'}$.
\par
(3.) This follows by repeated application of (2a) and (2b).
\end{proof}

For bigraded ideals in the polynomial ring
$R = \poly{K}{m}{n}$, we have:
\begin{pro}
\label{P:propweakreg}
Let $K$ be a field, let $I\subset R$ be any ideal
generated by bihomogeneous polynomials, let $\mc{I}$ be the
corresponding sheaf of ideals in $\mc{O}_{X}$.  The following
properties are equivalent.
\begin{enumerate}
\item The ideal $I$ is weakly $(p, p')$-regular in the sense of
Definition \ref{D:weakmodreg}.
\item The natural map $I_{p,p'} \to H^0(\mc{I}(p,p'))$ is
an isomorphism and $\mc{I}$ is $(p,p')$-regular in the sense of
Definition \ref{D:sheafreg}.
\item The natural map $I_{d,d'} \to H^0(\mc{I}(d,d'))$
is an isomorphism and $\mc{I}$ is
$(d,d')$-regular, for all $d \geq p, d' \geq p'$.
\end{enumerate}
\end{pro}
\begin{proof}
No loss in generality in assuming that $K$ is infinite, because
we may tensor the whole situation by the algebraic closure
of $K$. \newline
(I $\Rightarrow$ II) If $I$ is weakly $(p,p')$-regular in the sense
of Definition \ref{D:weakmodreg}, then we have
$\lclcoh{i}{\mbf{m}}{I}_{k,k'} = 0$ for  $(k,k')\in
Reg_{i-1}(p,p')$ for $i \geq 1$.
 But for an
ideal in a polynomial ring, we also have $\lclcoh{0}{\mbf{m}}{I} =
0$, since there are no $0$-divisors in ring $R$.  By Proposition
\ref{P:exact}, $H^i(\mc{I}(k,k'))=\lclcoh{i+1}{\mbf{m}}{I}_{k,k'} =
0$ for $i\ge 1$, $(k,k')\in Reg_{i}(p,p')$ for $i \geq 1$; and
$I_{p,p'}\cong H^0(\mc{I}(p,p'))$.
\par
(II $\Rightarrow$ I) If $\mc{I}$ is $(p,p')$-regular in the sense
of Definition \ref{D:sheafreg}, then $H^i(\mc{I}(k,k'))=0$ for
$(k,k')\in Reg_i(p,p')$, $i \geq 1$ by Proposition
\ref{P:greater}.   $I$ is an ideal, $\lclcoh{0}{\mbf{m}}{I} = 0$,
in particular, $H^0_{\mbf{m}}(I)_{k,k'}=0$ for $(k,k')\in
Reg '_{-1}(p,p')\cup Reg ''_{-1}(p,p')$. 
Since $I_{p,p'}\cong H^0(\mc{I}(p,p'))$, by
Proposition \ref{P:exact}, we have $H^1(\mc{I})_{p,p'}=0$, and
$H^{i+1}_{\mbf{m}}(I)_{k,k'}=H^i(\mc{I}(k,k'))=0$ for $(k,k')\in
Reg_i(p,p')$, $i \geq 1$.  Therefore $H^i_{\mbf{m}}(I)_{k,k'}=0$
for all $(k,k')\in Reg_{i-1}(p,p')$, $i  \geq 0$, i.e. $I$ is weak
$(p,p')$-regular in the sense of Definition \ref{D:weakmodreg}.
\par
(II $\Rightarrow$ III) follows from Proposition \ref{P:greater},
and Proposition \ref{P:span}.
\par
(III $\Rightarrow$ II) is obvious, we just take $d=p, d'=p'$.
\end{proof}

\section{Strong regularity for bigraded modules}
\label{S:modreg2}
From now on, $K$ is a field and
$R=K[x_0, \cdots, x_m, y_0, \cdots, y_n]
= K[x, y]$ is a polynomial
algebra, bigraded in the usual way. We will be
using the ideals
$(x)=(x_0,\cdots, x_m)$, $(y)=(y_0, \cdots, y_n)$,
$(x,y)=(x_0, \cdots, x_m,y_0, \cdots, y_n)$, and
$(xy)=\mbf{m}=(x_iy_j)$.
\par
In addition to the graded $K[x, y]$-module
$M^{\sharp}$ introduced above, we need to consider
graded modules as follows.
Fix $j'$, and let $M^{[1]}_{j'}=\oplus_j M_{j,j'}$,
which is a $K[x]=K[x_0,\cdots, x_m]$-module; fix $j$, and
let $M^{[2]}_{j}=\oplus_{j'} M_{j,j'}$,
which is a $K[y]=K[y_0, \cdots, y_m]$-module.
Observe that
\[
M = \bigoplus _{j' }M^{[1]}_{j'} =
\bigoplus _{j}M^{[2]}_{j}
\]
as $K[x]$-module (resp. as $K[y]$-module). Also,
each $\lclcoh{i}{(x)}{M^{[1]}_{j'}}$ is a graded
$K[x]$-module (resp. each $\lclcoh{i}{(y)}{M^{[2]}_{j}}$
is a graded $K[y]$-module), but both
$\lclcoh{i}{(x)}{M}$ and $\lclcoh{i}{(y)}{M}$ are bigraded
$K[x, y]$-modules.
\begin{align*}
\lclcoh{i}{(x)}{M} &= \bigoplus _{j'} \lclcoh{i}{(x)}{M^{[1]}_{j'}}\\
\lclcoh{i}{(y)}{M} &= \bigoplus _{j} \lclcoh{i}{(y)}{M^{[2]}_{j}}\\
\lclcoh{i}{(x)}{M}_{j, j'} &= \lclcoh{i}{(x)}{M^{[1]}_{j'}}_{j}\\
\lclcoh{i}{(y)}{M}_{j, j'} &= \lclcoh{i}{(y)}{M^{[2]}_{j}}_{j'}
\end{align*}

\begin{Def}\label{D:modreg}
Let $M$ be a bigraded $R$-module and let $d \geq 0$.
\begin{enumerate}
\item $M$ satisfies the \textit{vanishing condition $VC_d(p,p')$} if for
all $i \geq 0$
\begin{eqnarray*}
\lclcoh{i}{(x)}{M}_{k, k'}&= &H^i_{(x)}(M^{[1]}_{k'})_k=0,
\  \forall k \geq p+d+1-i, \forall k';\\
\lclcoh{i}{(y)}{M}_{k, k'}&= &H^i_{(y)}(M^{[2]}_k)_{k'}=0,
\  \forall k' \geq p'+d+1-i, \forall k;\\
& &H^i_{(x,y)}(M^{\sharp})_{k+k'}=0, \  \forall k+k' \geq
p+p'+d+1-i.
\end{eqnarray*}
\item $M$ is \textit{$(p,p')$-regular} if $M$ satisfies $VC_0(p,p')$.
\end{enumerate}
\end{Def}

\begin{rem}\label{R:modshiftreg}
For all $p,p'$, we have
\begin{itemize}
\item[1.]
If $M$ satisfies $VC_0(p,p')$, then $M$ satisfies $VC_d(p,p')$
for all $d\ge 0$.
\item[2.]
If $M$ satisfies $VC_d(p,p')$, then $M(a,b)$ satisfies
$VC_d(p-a,p'-b)$.
\item[3.]
For all $(\alpha, \alpha') \in DReg_d(p,p')$, if $M$ satisfies
$VC_0(\alpha, \alpha')$, then $M$ satisfies
$VC_d(p,p')$.
\end{itemize}
\end{rem}

\begin{pro}
\label{P:mzerozero} Let $R= K[x_0 , ..., x_m, y_0, ..., y_n]$ be a
bigraded polynomial algebra over a field $K$. Assume that $m, n
\ge 0$. Then $R$ is strongly $(0, 0)$-regular.
\end{pro}
\begin{proof}
If $R= K[z_1 , ..., z_s]$ is any polynomial algebra over a field
$K$, it is a classical fact, due essentially to Serre, that
$\lclcoh{i}{(z)}{R}_{k} = 0$ whenever $i+k \ge 1$. This verifies
the vanishing statement for $H^i_{(x,y)}(R^{\sharp})_{k+k'}$ for
$R = K[x, y]$. For the case $H^i_{(x)}(R^{[1]}_{k'})_{k}$, note
that
\[
R^{[1]}_{k'} = \bigoplus _{| \beta | = k'}K[x]y^{\beta}.
\]
As each term is a free module over $K[x]$ and local cohomology
commutes with direct sum, the requisite vanishing follows from
Serre's result.
\end{proof}

\begin{pro}\label{P:Mvanreg}
If a bigraded $R$-module $M$ satisfies $VC_d(p,p')$, then
$H^i_{(xy)}(M)_{k,k'}=0$ for all $(k,k') \in
Reg_{0}(p+d+1-i,p'+d+1-i), 0 \leq i \leq d+2$;
and for all $(k,k')\in Reg_{i-d-1}(p,p'), i \geq d+2$.
\end{pro}

\begin{proof}
By the Mayer-Vietoris sequence, we have
\[
\begin{CD}
 \lclcoh{i}{(x, y)}{M}
@>>> \lclcoh{i}{(x)}{M} \oplus
\lclcoh{i}{(y)}{M} @>>> \lclcoh{i}{(xy)}{M}@>>>
\lclcoh{i+1}{(x, y)}{M}
\end{CD}
\]
(see \cite[Exercise 2.4, Ch. III, p. 212]{H};
Note that if $Y_1 = V(x)$, $Y_2 = V(y)$, then
$Y_1 \cup Y_2 = V(xy)$ and $Y_1 \cap Y_2 = V(x, y)$ as subsets
of $C = \text{Spec}(K[x, y])$. )
Assuming that $M$ satisfies $VC_d(p,p')$ we see that
$H^i_{(xy)}(M)_{k,k'}=0$ for all $(k,k')$ that satisfy the
inequalities:

\[
k \ge p+d+1-i, \ \ k' \ge p'+d+1-i,\
k+k' \ge p+p'+d-i.
\]
If $0 \le i \le d+2$ the last condition above is redundant, and
so we obtain vanishing in the region described by the
first two inequalities, which is just $ Reg_{0}(p+d+1-i,p'+d+1-i)$.
If  $i \geq d+2$, these three inequalities describe
$Reg_{i-d-1}(p,p')$.
\end{proof}

\begin{cor}
\label{P:strweak}
If $M$ is strongly $(p, p')$-regular, then it is weakly $(p, p')$-regular.
\end{cor}
\begin{proof}
We have $\lclcoh{i}{\mbf{ m}}{M}_{k. k'} = 0$ for all
$(k, k') \in Reg_{i-1}(p, p')$, according to the proposition,
whenever $i\ge 2$. For $i = 0, 1$, this is zero for
$(k, k') \in Reg_{0}(p+1-i, p'+1-i)$, but these are exactly
the regions $Reg_{-1}(p, p')$ and  $Reg_{0}(p, p')$. Thus
we have the conditions for weak $(p, p')$-regularity.
\end{proof}

\begin{rem}
\label{R:shift}
 If $M$ is strongly (resp. weakly) $(p, p')$-regular, then
 $M(a, b)$ is strongly (resp. weakly) $(p-a, p'-b)$-regular.
\end{rem}

\begin{pro} \label{P:gen}
If a finitely generated bigraded $R$-module
$M$ satisfies $VC_d(p,p')$, then $M$ is
generated by elements of bidegree $(k, k') \in DReg_d(p,p')$.
\end{pro}

\begin{proof}
Let $A$ be a homogeneous algebra in the sense of
Ooishi's paper (see introduction to \cite{O}), with
maximal ideal $P$.
If $N$ is a finitely generated graded module over $A$,
then \cite[Thm. 2]{O} asserts that
if $\lclcoh{i}{P}{N}_k = 0$ for all $i+k \ge m+1$, then
$N$ is generated in degrees $\le m$.
\par
We first apply this to the graded module $N = M^{\sharp}$
over the graded ring $A = R^{\sharp}$.
Since $M$ satisfies $VC_d(p,p')$, we have
\[
 H^i_{(x,y)}(M^{\sharp})_{k+k'}=0, \  \forall k+k' \geq
        p+p'+d+1-i
\]
so that by the previous remark, $M^{\sharp} $ can be
generated by elements of degree $\le p+p'+d$.
This means that the bigraded $M$ can be generated
by bihomogeneous elements of bidegree $(k, k')$
with $k + k' \le p+p'+d$. Now let $A = K[x]$, and for a
fixed $k'$, regard $N = M^{[1]}_{k'}$ as an $A$-module.
That  $M$ satisfies $VC_d(p,p')$ means here that
\[
 H^i_{(x)}( M^{[1]}_{k'})_k=0, \  \forall k \geq p+d+1-i
\]
and thus by Ooishi's result, that $M^{[1]}_{k'}$ can
be generated as $ K[x]$-module by elements of degree
$\le p+d$. This being true for every $k'$, we see that
\[
R_{s, 0}M_{p+d, k'} = M_{p+d+s, k'}\quad \text{for all }
s \ge 0, k'.
\]
Similar reasoning applied to  $M^{[2]}_{k}$ as an $K[y]$-module
leads to
\[
R_{0, s}M_{k, p'+d} = M_{k, p'+d+s}\quad \text{for all }
s \ge 0, k.
\]
Combining this information gives that $M$ can be
generated by bihomogeneous elements of degree
$(k, k')$ where
\[
k \le p+d, \quad k' \le p'+d\quad k+k' \le p+p'+d
\]
This is the description of the region
$DReg_{d}(p, p')$.
\end{proof}

If $M_d$ is a bigraded $R$ module which satisfies $VC_d(p,p')$, by
Proposition \ref{P:gen}, $M_d$ is generated by elements of
bidegree $e_{\alpha, d}=(\alpha_d, \alpha'_d) \in DReg_d(p,p')$.
We can find an exact sequence:
\[\begin{CD}  0@>>>M_{d+1} @>>>
\displaystyle{\bigoplus_{\alpha=1}^{r_d}} Re_{\alpha, d}
@>{\phi_d}>> M_d @>>> 0,\end{CD}\] where $M_{d+1}=\ker \phi_d$.

\begin{pro}\label{P:Induct}
Let $M_d$ be as above.  If $M_d$ satisfies $VC_d(p,p')$, then $M_{d+1}$
satisfies $VC_{d+1}(p,p')$, and therefore are generated by elements of
bidegree in $DReg_{d+1}(p,p')$.
\end{pro}

\begin{proof}

For the case $i=0$, we have an injection
\[
\lclcoh{0}{(x)}{M_{d+1}} \subset
\displaystyle{\bigoplus_{\alpha=1}^{r_d}}
H^{0}_{(x)}
(R)_{k-\alpha _{d},k'-\alpha'_{d}} = 0,
\]
so we can assume that $i \ge 1$.
Consider the local cohomology sequence
with $I = (x)$ of the above exact
sequence:
\[
\begin{CD}
H^{i-1}_{I}(M_d)_{k,k'}
  @>>> H^{i}_{I}(M_{d+1})_{k,k'} @>>>
\displaystyle{\bigoplus_{\alpha=1}^{r_d}}
H^{i}_{I}
(R)_{k-\alpha _{d},k'-\alpha'_{d}}
\end{CD}
\]
Suppose that $k+i \ge p+(d+1)+1$. Then the left-hand side above
vanishes by assumption on $M_d$, because
 $k+(i-1) \ge p+d+1$. That $(\alpha _{d}, \alpha'_{d})$
belongs to $DReg_d(p, p')$ means that
$\alpha _d \le p+d$. Thus,
$k -\alpha _d + i \ge 2 $, and since $R$ is $(0, 0)$-regular
by Proposition (\ref{P:mzerozero}),  the last term vanishes.
\par
By similar reasoning, we get the vanishing of
$\lclcoh{i}{(y)}{M_{d+1}}_{k, k'}$  for $k'+i \geq
p'+(d+1)+1$, for all $k$.
\par

Now look at the local cohomology sequence with
$I = (x, y)$. Again we may assume that $i \ge 1$.
If $(k, k')$ satisfies $k+k'+i \geq p+p'+(d+1)+1$,
our assumption on $M_d$ shows the vanishing of the
left-hand side because $k+k'+(i-1) \geq p+p'+d+1$.
 That $(\alpha _{d}, \alpha'_{d})$
belongs to $DReg_d(p, p')$ means that
$\alpha _d + \alpha'_d \le p+p' +d$, so that
 $k+k'-\alpha_d-\alpha'_d+ i \geq 2$. Thus the right-hand side
vanishes because $R$ is $(0, 0)$-regular.
\par
In all three cases we have verified vanishing
in the appropriate region to satisfy
$VC_{d+1} (p, p')$.
\end{proof}
Conversely:
\begin{pro}\label{P:Inverseinduct}
Let $M_{d+1}$ be a finitely generated bigraded $R$-module.
If $M_{d+1}$ satisfies
$VC_{d+1}(p,p')$ and if there is an exact sequence:
\[
\begin{CD}  0@>>>M_{d+1} @>>>
\displaystyle{\bigoplus_{\alpha=1}^{r_d}} Re_{\alpha, d}
@>{\phi_d}>> M_d @>>> 0,
\end{CD}
\]
where $M_{d+1}=\ker \phi_d$,
and $e_{\alpha,d}
= (\alpha _d , \alpha'_d)\in DReg_d(p,p')$, then $M_{d}$ satisfies
$VC_{d}(p,p')$.  Therefore $M_d$ is generated by elements of
bidegree in $DReg_{d}(p,p')$.
\end{pro}

\begin{proof}
Let $I$ be any one of the ideals $(x)$, $(y)$, $(x, y)$.
Look at the segment of the local cohomology
sequence associated with the above exact sequence:
\[
\begin{CD}
\displaystyle{\bigoplus_{\alpha=1}^{r_d}}
H^{i}_{I}
(R)_{k-\alpha _{d},k'-\alpha'_{d}}@>>>
H^{i}_{I}(M_d)_{k,k'}
  @>>> H^{i+1}_{I}(M_{d+1})_{k,k'}
\end{CD}
\]
Let $I = (x)$.
and suppose $k+i \geq p+d+1$,  we have
$k+(i+1) \geq p+(d+1)+1$,
and the last group vanishes by assumption on $M_{d+1}$.
Also, in this region, $k-\alpha_d+i \ge 1$, and the
first term vanishes by Proposition (\ref{P:mzerozero}).

By similar reasoning, we obtain the vanishing of
$H^{i}_{(y)}(M_d)_{k,k'}$ if  $k'+i \geq p'+d+1$.
\par
For $I= (x, y)$, suppose $k+k'+i \geq p+p'+d+1$.  We have
$k+k'+(i+1) \geq p+p'+(d+1)+1$, so that the last group
vanishes by assumption on $M_{d+1}$.
Also, because
$\alpha _d + \alpha ' _d \le p+p' +d$, we get
$k+k'-\alpha_d-\alpha'_d -i\geq 1$, last term vanishes
because $R$ is $(0, 0)$-regular.
\par
In all three cases we have verified vanishing
in the appropriate region to satisfy
$VC_d (p, p')$.
\end{proof}

We  prove some equivalent conditions for regularity of a module,
as in \cite{BS} and \cite{BM}. In the formulation below, $R$ is a
polynomial algebra over $K$ in two sets of variables $x$ and $y$
bigraded in the usual way. We assume both variable sets are
nonempty.

\begin{thm}\label{T:BayerMum} Let $M$
be a finitely generated bigraded module over $R$.  The
following properties are equivalent.
\begin{enumerate}
\item $M$ is $(p, p')$-regular in the sense of
definition \ref{D:modreg}.
\item The  minimal resolution of $M$ by free bigraded
$R = K[x,y]$-modules:
\[
\begin{CD}
0 @>>> \displaystyle{\bigoplus_{\alpha=1}^{r_s}}
R e_{\alpha, s} @>>> \cdots @>>>
 \displaystyle{\bigoplus_{\alpha=1}^{r_0}}R e_{\alpha, 0} @>>> M @>>> 0,
\end{CD}
\]
satisfies $e_{\alpha, d}=(\alpha_d, \alpha'_d) \in DReg_d(p,p').$
\item There exists a free resolution with the properties above.
\end{enumerate}
\end{thm}

\begin{proof}
(I $\Rightarrow$ II) Let $M_0 = M$. We will inductively construct
a sequence of bigraded modules $M_d$ that satisfy
$VC_d(p,p')$  and that fit into an exact sequence
\begin{equation}
\label{E:induct}
\begin{CD}
0@>>>M_{d+1} @>>> \displaystyle{\bigoplus_{\alpha=1}^{r_d}} Re_{\alpha, d}
@>{\phi_d}>> M_d @>>> 0
\end{CD}
\end{equation}
where  $e_{\alpha, i}=(\alpha_i, \alpha'_i) \in DReg_d(p,p')$.
By Proposition (\ref{P:Induct}),
we know that $M_{d+1}$ will satisfy $VC_{d+1}(p,p')$ and therefore
we can find generators for it whose  bidegrees are in
$DReg_{d+1}(p,p')$. In other words, we may construct the
above exact sequence with but $d$ replaced by $d+1$. By Hilbert's
syzygy theorem, $M_d$ will become a free bigraded module, with
generators in $DReg_{d}(p,p')$, and by splicing these short
sequences together, we get our resolution. We can start this induction
at $d=0$, because by hypothesis, $M=M_0$ is $(p,p')$-regular,
and by Proposition (\ref{P:gen}), we know
$M_0$ is generated by elements whose  bidegrees are in $DReg_0(p,p')$.
\par
(II $\Rightarrow$ III) is trivial.
\par
(III $\Rightarrow$ I) Break the given resolution into short
sequences as in equation (\ref{E:induct}) above. We will
show by descending induction on $d$ that $M_d$ satisfies
$VC_d(p,p')$. Since the last stage of this, namely $M_0$,
is the module $M$ itself, we will be done, since the condition
$VC_0(p,p')$ is exactly $(p, p')$-regularity. The starting
point of the induction is the extreme left-hand term of
the resolution
$M_s=\oplus_{\alpha=1}^{r_s}Re_{\alpha_s, \alpha'_s}$. Because
$R$ is $(0, 0)$-regular by Proposition (\ref{P:mzerozero}), and
because of Remark (\ref{R:modshiftreg}), we see that
$M_s$ satisfies $VC_s(p,p')$. If $d < s$ and we assume by induction
that $M_{d+1}$ satisfies $VC_{d+1}(p,p')$, from the exact sequence
(\ref{E:induct}) and Proposition
(\ref{P:Inverseinduct}), we find that $M_d$ satisfies $VC_d(p,p')$,
verifying the induction step.
\end{proof}

\begin{cor}
\label{C:BayerMum} Any finitely generated bigraded module
over $K[x, y]$ is $(p, p')$-regular for some $p, p'$.
\begin{proof}
Look at the minimal free bigraded resolution of $M$, which
we know exists and is unique up to isomorphism. Whatever are
the bidegrees $e_{\alpha, d}$ of the generators of the
various terms in this, it is rather clear that by taking
$p$ and $p'$ sufficiently large, for all $d$ these will
belong to the region $DReg_d(p, p')$.
\end{proof}
\end{cor}

\begin{rem} \label{R:IReduced}
Let $I \subset K[x_0, y_0, \cdots, y_n]$ be an ideal such that
$I=x_0^m J$ where $J \subset K[y_0, \cdots, y_n]$ is a homogeneous
ideal, where $K$ is an infinite filed.  Then  $I$ is $(p,p')$-regular if and only if $p \ge m$ and
$J$ is $p'$-regular.
\end{rem}

\begin{proof}
Suppose $J$ is $p'$-regular and $p \ge m$, we would like to show
that $I$ is $(p,p')$-regular. Let $R=K[y_0, \cdots, y_n]$, and
take a minimal free resolution of $J$ as follows:
\begin{equation}
\label{E:MinJ}
 \begin{CD}
0 @>>>  \displaystyle{\bigoplus_{\alpha=1}^{r_s}} R e_{\alpha, s}
@>>> \cdots @ > d_1 >>
 \displaystyle{\bigoplus_{\alpha=1}^{r_0}}
R e_{\alpha,0} @> d_0
>>  R @>>>  R/J  @>>>  0
\end{CD}.
\end{equation}
 Since $J$ is $p'$-regular, then we have $e_{\alpha,i} \le p'+i.$
Also, note the map $d_i$ is represented by a matrix.  Since the
free resolution is minimal, then the matrix has no entry in $K^*$
where $K^*=K \setminus \{0\}$. \cite[Proposition 11.5]{SR}

We can break this exact sequence into two exact sequences:
\begin{equation} \label{E:kerJ}
\begin{CD}
 0 @>>> C_n
@>>> \cdots @>>> C_1 @>>> \ker(d_0)
@>>> 0,
\end{CD}
\end{equation}

 and
\begin{equation}
\label{E:J}
\begin{CD}
0
@>>> \ker(d_0) @>>>  \displaystyle{\bigoplus_{\alpha=1}^{r_0}}R
e_{\alpha, 0} @> d_0
>>  R @>>>   R/J  @>>>  0
\end{CD}.
\end{equation}

If we let $S=K[x_0, y_0, \cdots, y_n]$, we have
\begin{equation}
\label{E:I}
\begin{CD}
0 @>>> \ker(x^m_0 d_0) @>>>
\displaystyle{\bigoplus_{\alpha=1}^{r_0}}S
e'_{\alpha, 0} @> d_0
>>  S @>>>   S/I  @>>>  0
\end{CD}.
\end{equation}

If we tensor the exact sequence \eqref{E:kerJ} with $S$ over $R$,
since $S$ is flat, then we will have an exact sequence:
\begin{equation} \label{E:TensorJ}
\begin{CD}
 0  @>>> C_n \otimes S
@>>> \cdots @>>> C_1 \otimes S @>>> \ker(d_0)
\otimes S @>>> 0.\end{CD}
\end{equation}
Note, at each stage, the matrix which represents the map has no
entry in $K^*$. Since $\ker(x^m_0 d_0)=S \otimes_{R} \ker(d_0)$,
we can piece exact sequence \eqref{E:I} and \eqref{E:TensorJ}
together, we will form a free resolution of $I$ as follow:
\begin{equation}
\label{E:TensorI}
 \begin{CD}
0  @>>> C_n \otimes S
@>>> \cdots @>>> C_1 \otimes S @>>>
 \displaystyle{\bigoplus_{\alpha=1}^{r_0}}S e'_{\alpha, 0} @> d_0 >>  S @>>>
S/I  @>>> 0
\end{CD}.
\end{equation}
This free resolution is minimal, since the matrix represents the
map has no entry in $K^*$.  And we can rewrite the minimal free
resolution \eqref{E:TensorI} as follow:
\begin{equation}
\label{E:MinI}
\begin{CD}
 0 @>>>  \displaystyle{\bigoplus_{\alpha=1}^{r_s}} S e'_{\alpha, s}
@>>> \cdots @> d'_1 >>  \displaystyle{\bigoplus_{\alpha=1}^{r_0}}
S e'_{\alpha,0} @> d'_0 >> S @>>> S/I @>>> 0
\end{CD}
\end{equation}
where $d'_0=x^m_0 d_0$, and $d_1'=d_1$, and $e'_{\alpha, i} =(-m,
e_{\alpha, i})$.  If $m\leq p$ and $e_{\alpha, i} \leq p'+i$, by
the equivalent relation of minimal free resolution and
$(p,p')$-regular, we know that $I$ is $(p, p')$-regular.

\par
On the other hand, suppose $I$ is $(p,p')$-regular,  there is a
minimal free resolution of $I$ as \eqref{E:MinI}, where
$e'_{\alpha, i}=(a_{\alpha, i}, e_{\alpha, i})$ and $a_{\alpha, i}
\leq p$ and $e_{\alpha, i} \le p'+i$. Note, at each stage, the
matrix represents the map has no  entry in $K^*$.  And we can
split the free resolution into two exact sequences: the free
resolution of $I$ \eqref{E:I} and

\begin{equation}\label{E:KerI1}
\begin{CD}
0  @>>> D_n @ > d'_n >> \cdots @>>> D_1 @ > d'_1 >>
\ker(x^m_0 d_0) @>>> 0
\end{CD},
\end{equation}

We always have a resolution of $J$ as \eqref{E:J}. Since
$\ker(x^m_0 d_0)=S \otimes_{R} \ker(k_0)$, we will have an exact
sequence as follow:
\begin{equation}
\label{E:TensorI1}
\begin{CD}
0  @>>> C_n @ > d_n >> \cdots @>>> C_1 @ > d_1 >>
\ker(d_0) @>>> 0
\end{CD},
\end{equation}
where $d_i=d'_i$. We can piece the two exact sequences
\eqref{E:TensorI1} and \eqref{E:J} together to get:
\[
\begin{CD}
0  @>>> C_n @ > d_n >> \cdots @>>> C_1 @ > d_1 >>
 \displaystyle{\bigoplus_{\alpha=1}^{r_0}}R e_{\alpha, 0} @> d_0
>>  R @>>>   R/J  @>>>  0
\end{CD},
\]
which can be written as \eqref{E:MinJ}. Since the matrix
represents $d_i$ has no  entry in $K^*$, the free resolution
\eqref{E:MinJ} is minimal, and $e_{\alpha,i} \le p'+i$. As is
well-known, the existence of a free resolution of this type
implies that $J$ is $p'$-regular.
\end{proof}

\section{Resolutions}
\label{S:resolutions} Let $R = K[x_0, ..., x_m, y_0, ..., y_n]$, a
polynomial algebra over a field, bigraded in the usual way. Let
$\mbf{m} = (x_i y_j)$, the irrelevant ideal. We will define a
complex that allows the computation of the local cohomology
modules $\lclcoh{i}{\mbf{m}}{M}$. Recall that for any ideal $I$ in
a ring $R$ we have
\[
\lclcoh{i}{I}{M} = \injlim _{\nu} \text{Ext}^i _{R}(R /I^{(\nu)},
M)
\]
where $I^{(\nu)}$ is any sequence of ideals cofinal with the
collection of powers $I^{\nu}$. If $R = K[z_1, ..., z_n]$ is a
polynomial ring and $I = (z_1, ..., z_n)$, then we can take
$I^{(\nu)}= (z_1^{\nu}, ..., z_n^{\nu}) $. Since this is generated
by a regular sequence, we can compute the $\text{Ext}$ groups by
using the Koszul complex on $z^{\nu} = \{ z_1^{\nu}, ...,
z_n^{\nu}\}$ as a free resolution of $R /I^{(\nu)}$:
\[
K_{\ast}(z^{\nu}) \longrightarrow R /I^{(\nu)}.
\]
This means that we have a free resolution by the truncated
complex:
\[
K_{\ge 1}(z^{\nu}) \longrightarrow I^{(\nu)}
\]
Let, as before, $\mbf{m}_1 = (x_0, ..., x_m)\subset R_1 = K[x_0,
..., x_m]$, $\mbf{m}_2 = (y_0, ..., y_n)\subset R_2 =  K[x_0, ...,
y_n]$, with  $\mbf{m}_1^{(\nu)} = (x_0^{\nu}, ..., x_m^{\nu})$,
with similar notation for $\mbf{m}_2^{(\nu)}$, $\mbf{m}^{(\nu)}$
If $M_1$ is an $R_1$-module and $M_2$ is an $R_2$-module, we let
$M_1 \boxtimes M_2$ be the $R$-module $(M_1 \otimes _{R_1}
R)\otimes _{R}(M_2 \otimes _{R_1} R)$.

\begin{lem}
\label{L:m1m2} There is a canonical isomorphism $\mbf{m}_1^{(\nu)}
\boxtimes \mbf{m}_2^{(\nu)} \cong \mbf{m}^{(\nu)}$ for every $\nu
\ge 0$.
\end{lem}
\begin{proof}
Clearly there is an epimorphism $\mbf{m}_1^{(\nu)} \boxtimes
\mbf{m}_2^{(\nu)} \to \mbf{m}^{(\nu)}$, so the issue is the
injectivity of this map. Because $R$ is a flat $R_1$-module, we
have an exact sequence
\[
\begin{CD}
0 @>>> \mbf{m}_1^{(\nu)} \otimes _{R_1}R=
 \mbf{m}_1^{(\nu)}R
@>>> R @>>> R/ \mbf{m}_1^{(\nu)}R @>>> 0.
\end{CD}
\]
Applying $- \otimes _{R} \mbf{m}_2^{(\nu)}R $, we see that the
lemma holds if we can show that
\[
\text{Tor}_{1}^{R}( \mbf{m}_2^{(\nu)}R, \, R/ \mbf{m}_1^{(\nu)}R )
= 0
\]
We can use the truncated Koszul complex on  $\mbf{m}_2^{(\nu)}R$
to compute the Tor. Since the elements $y_0 ^{\nu}, ..., y_n
^{\nu}$ are a regular sequence on $R/ \mbf{m}_1^{(\nu)}R$, this
Tor vanishes as claimed.
\end{proof}
As we have mentioned, there are Koszul resolutions
\[
K_{\ge 1}(x^{\nu}) \longrightarrow \mbf{m}_1 ^{(\nu)}, \quad K_{\ge
1}(y^{\nu}) \longrightarrow \mbf{m}_2 ^{(\nu)}
\]
\begin{pro}
\label{P:resolution}
\[
K_{\ge 1}(x^{\nu})\boxtimes K_{\ge 1}(y^{\nu})\longrightarrow
  \mbf{m}_1 ^{(\nu)}\boxtimes\mbf{m}_2 ^{(\nu)} =  \mbf{m} ^{(\nu)}
\]
is a finite resolution by free $R$-modules of finite type.
\end{pro}
\begin{proof}
The only thing to be proved is that it actually is a resolution.
Since $R$ is a flat $R_j$-module, $j=1, 2$, we get free
resolutions of $R$-modules:
\[
A_{\ast}: = K_{\ge 1}(x^{\nu})\otimes _{R_{1}}R \longrightarrow
\mbf{m}_1 ^{(\nu)}\otimes _{R_{1}}R =\mbf{m}_1 ^{(\nu)}R , \quad
B_{\ast}: = K_{\ge 1}(y^{(\nu)}) \otimes _{R_{2}}R \longrightarrow
\mbf{m}_2 ^{(\nu)}\otimes _{R_{2}}R = \mbf{m}_2 ^{(\nu)}R
\]
Clearly we have an epimorphism $K_{\ge 1}(x^{\nu})\boxtimes K_{\ge
1}(y^{\nu}) = A_{\ast}\otimes _{R}B_{\ast} \to \mbf{m}^{(\nu)}$, so
we must prove that
\[
\text{H}_i (A_{\ast}\otimes _{R}B_{\ast}) = 0, \text{ for all } i
\ge 1.
\]
Since the $A_r$ are free $R$-modules, we have the K\"unneth
spectral sequence (see \cite[Ch. XVII, 5.2a]{CE},
\cite[I, 5.5.2]{rG}, \cite[Lemme (1.1.4.1)]{aGIII}):
\[
\text{E}_{p, q}^2 = \text{H}_p (A_{\ast}\otimes _{R}\text{H}_q
(B_{\ast}))) \Longrightarrow\text{H}_{p + q} (A_{\ast}\otimes
_{R}B_{\ast})
\]
These are zero when $q \ge 1$ since $B_{\ast} $ is a resolution.
We must see that they are $0$ when $p \ge 1$, or that
\[
\text{H}_p (A_{\ast}\otimes _{R} \mbf{m}_2 ^{(\nu)}R  )=
\text{H}_{p+1} (K_{\ast}(x^{\nu})\otimes _{R} \mbf{m}_2 ^{(\nu)}R  )
= \text{Tor}_{p+1}^{R}( R/ \mbf{m}_1^{(\nu)}R,\,
  \mbf{m}_2^{(\nu)}R ) = 0,
\ p \ge 1.
\]
This follows because as noted in the previous lemma, the elements
$y_0 ^{\nu}, ...,y_n ^{\nu} $ form a regular sequence in $ R/
\mbf{m}_1^{(\nu)}R$.
\end{proof}


\begin{thebibliography}{30}
\bibitem{BM} D. Bayer and D. Mumford, {\em What can be computed in
algebraic geometry}, Computational Algebraic Geometry and
Commutative Algebra.  Cambridge Univ. Press, Cambridge, 1993,
1-48.

\bibitem{BS} D. Bayer and M. Stillman, {\em A criterion for
detecting m-regularity}, Invent. Math. \textbf{87} (1987), 1-11.

\bibitem{CE} H. Cartan and S. Eilenberg, {\em Homological Algebra},
   Princeton University Press., Princeton, N.J. ,  1956.
  
\bibitem{EG} D. Eisenbud and S. Goto, {\em Linear free resolution and minimal multiplicity},
   J. Algebra, \text{88} (1984), 89-133.
   
\bibitem{rG} R. Godement, {\em Topologie alg\'ebrique et th\'eorie des faisceaux},
 Actualit{\'e}s Sci. Ind. No. 1252. Publ. Math. Univ. Strasbourg. No. 13, Hermann, Paris, 1958.

\bibitem{GW1} S. Goto and K. Watanabe, {\em On graded rings. {I}.}, J. Math. Soc. Japan,
 \text{30} (1978), 179-213.
 
 \bibitem{GW2} S. Goto and K. Watanabe, {\em On graded rings. {II. ($Z\sp{n}$-graded rings)}},
  Tokyo J. Math., \text{1} (1978), 237-261.
  
\bibitem{aGIII} A. Grothendieck, {\em {\'El\'ements de g\'eom\'etrie alg\`ebrique. III. 
               \'Etude cohomologique des faisceaux coh\'erents. I.}}, Inst. Hautes {\'E}tudes Sci. Publ. Math.,
 \text{11} (1961).
 
\bibitem{aG} A. Grothendieck, {\em Cohomologie locale des faisceaux cohérents et
                  th\'eor\`emes
                  de {L}efschetz locaux et globaux ({SGA2}).}, North-Holland Publishing Co.,
                  Amsterdam; Masson \& Cie, \'Editeur, Paris, 1968.

\bibitem{H} R. Hartshorne,  {\em Algebraic Geometry},
 Graduate Texts in Mathematics,
\textbf{52}, Springer-Verlag, 1977.


\bibitem{HY} E. Hyry, {\em The diagonal subring and the
Cohen-Macauly property of a multigraded ring}, Trans. Amer. Math.
Soc. \textbf{351}, no 6, 2213-2232.


\bibitem{JK}  B. Johnston and D. Katz, {\em Castelnuovo regularity and graded rings associated to an ideal},
 Proceedings of the American Mathematical Society., \text{123} no 3, (1995), 727-734.
 
\bibitem{DM} D. Mumford, {\em Lectures on curves on an algebraic surface},
Princeton University Press, Princeton, New Jersey, 1966.

\bibitem{O} A. Ooishi, {\em Castelnuovo's regularity of graded rings and modules},
Hiroshima Math. J. , \text{12} (1982), 627-644.

\bibitem{QB} P. Quartararo and H. S. Butts, {\em Finite unions of ideals and modules. },
Proc. Amer. Math. Soc.,  \text{52} (1975), 91-96.

\bibitem{SW} J. H. Sampson and G. Washnitzer,
{\em A {K\"u}nneth formula for coherent algebraic sheaves},
Illinois J. Math., \textbf{3}, 1959, 389-402.

\bibitem{Serre1} J.-P. Serre, {\em {Faisceaux alg\`ebriques coh\'erents}},
Ann. of Math., \text{61} (1955), 197-278.
  
 \bibitem{SR} R. P. Stanley, {\em Combinatorics and commutative algebra},
Birkh$\ddot{a}$user, 1983.
 
  
\end{thebibliography}
\end{document}